\documentclass[11pt]{article}

	\usepackage{amsmath, amssymb, amscd, latexsym, epsfig}
\newcommand{\arrow}{\rightarrow}

\newcommand{\bb}{\mathbb}
\newcommand{\cx}{{\bb C}}

\newcommand{\half}{{\bb H}}
\newcommand{\integers}{{\bb Z}}

\newcommand{\ratls}{{\bb Q}}
\newcommand{\reals}{{\bb R}}
\newcommand{\proj}{{\bb P}}

\thicklines
\renewcommand{\bold}[1]{\smallskip \noindent {\bf \boldmath #1 }\nopagebreak[4]}

\newcommand{\qed}{\nopagebreak[4]\hfill
\rule{2mm}{2.5mm} \bigskip \pagebreak[2]}

\newcommand{\Arrow}{\longrightarrow}
\newcommand{\brackets}[1]{\langle #1 \rangle}

\newcommand{\compos}{\circ}

\newcommand{\isom}{\cong}

\newcommand{\mem}{\in}

\newcommand{\notmem}{\not\in}

\newcommand{\AND}{\;\;\;\text{and}\;\;\;}

\newcommand{\st}{\: : \:}         %Such that

\newcommand{\ratlsbar}{{\overline{\ratls}}}

\newcommand{\chat}{{\widehat{\cx}}}

\newcommand{\Aut}{\operatorname{Aut}}

\newcommand{\id}{\operatorname{id}}

\newcommand{\mult}{\operatorname{mult}}

\newcommand{\zed}{\integers}

\newcommand{\cP}{{\cal P}}
\newcommand{\cQ}{{\cal Q}}

\newcommand{\cT}{{\cal T}}

\newtheorem{theorem}{Theorem}[section]
\newtheorem{prop}[theorem]{Proposition}
\newtheorem{lemma}[theorem]{Lemma}
\newtheorem{cor}[theorem]{Corollary}

\newtheorem{quest}[theorem]{Question}

\makeatletter
\def\cleardoublepage{\clearpage\if@twoside \ifodd\c@page\else
    \thispagestyle{plain}\hbox{}\newpage\if@twocolumn\hbox{}\newpage\fi\fi\fi}

\def\ps@headings{\let\@mkboth\markboth
  \def\@oddfoot{}%
  \def\@evenfoot{}%
  \def\@evenhead{\small \sc\thepage\hfil\leftmark}%      Left heading.
  \def\@oddhead{\small \sc \rightmark\hfil\thepage}%Right heading
  \def\chaptermark##1{{%  Expand the \ifnum, not the ##1
    \edef\@tempa{\ifnum \c@secnumdepth >\m@ne \@chapapp\ \thechapter. \fi}%
    \expandafter \markboth \expandafter{\@tempa ##1}{}}}%
  \def\schaptermark##1{\markboth {##1}{##1}}%
  \def\sectionmark##1{{%  Expand the \ifnum, not the ##1
    \edef\@tempa{\ifnum \c@secnumdepth >\z@ \thesection. \fi}%
    \expandafter \markright \expandafter{\@tempa ##1}}}}

\def\thebibliography#1{\section*{References\@mkboth
 {References}{References}}\list
 {[\arabic{enumi}]}{\settowidth\labelwidth{[#1]}\leftmargin\labelwidth
 \advance\leftmargin\labelsep
 \usecounter{enumi}}
 \def\newblock{\hskip .11em plus .33em minus .07em}
 \sloppy\clubpenalty4000\widowpenalty4000
 \sfcode`\.=1000\relax}

\newif\if@restonecol
\def\theindex{\@restonecoltrue\if@twocolumn\@restonecolfalse\fi
\columnseprule \z@
\columnsep 35pt\twocolumn[\@makeschapterhead{Index}]
 \@mkboth{Index}{Index}\thispagestyle{plain}\parindent\z@
 \parskip\z@ plus .3pt\relax\let\item\@idxitem}
\def\@idxitem{\par\hangindent 40pt}

\def\endtheindex{\if@restonecol\onecolumn\else\clearpage\fi}

\def\footnoterule{\kern-3\p@ 
 \hrule width .4\columnwidth 
 \kern 2.6\p@} 
\@addtoreset{footnote}{chapter} 
\long\def\@makefntext#1{\parindent 1em\noindent 
 \hbox to 1.8em{\hss$^{\@thefnmark}$}#1}

\@addtoreset{equation}{section}

\renewcommand{\l@section}{\@dottedtocline{0}{1.5em}{2.3em}}
\renewcommand{\l@subsection}{\@dottedtocline{1}{3.8em}{3.2em}}
\renewcommand{\l@subsubsection}{\@dottedtocline{2}{7.0em}{4.1em}}

\makeatother

\usepackage{caption,xy}
\xyoption{all}

\title{\vspace{-1in}
	{\bf On the postcritical set of a rational map}
\vspace{.2in}}

\author{Laura G. DeMarco, Sarah C. Koch and Curtis T. McMullen}

\date{20 September 2017}

\begin{document}
\maketitle

\begin{abstract} 
The postcritical set $P(f)$ of a rational map $f:\proj^1\to \proj^1$ is the smallest 
forward invariant subset of $\proj^1$ that contains the critical values of $f$.
In this paper we show that every finite set $X\subset \proj^1(\ratlsbar)$ can 
be realized as the postcritical set of a rational map.
We also show that every map $F:X\to X$ defined on a finite set  
$X\subset \proj^1(\cx)$ can be realized by a rational map $f:P(f)\to P(f)$, 
provided we allow small perturbations of the set $X$.
The proofs  involve Belyi's theorem and iteration on Teichm\"uller space. 
\end{abstract}

\setcounter{tocdepth}{1}
\tableofcontents

\vfill \footnoterule \smallskip
{\footnotesize \noindent
 Research supported in part by the NSF.
 %Typeset \today.
 }

\thispagestyle{empty}
\setcounter{page}{0}

\newpage

\section{Introduction}

Let $f:\proj^1\to \proj^1$ be a rational map on the Riemann sphere $\proj^1 = \proj^1(\cx)$ 
of degree $d\geq 2$.
Let $C(f) \subset \proj^1$ denote the set of critical points of $f$,
and let $V(f)$ denote the set of critical values.
The {\em postcritical set} of $f$ is defined by
$$P(f) = \bigcup_{n\geq 0} f^n(V(f)).$$

A rational map is {\em postcritically finite} if $|P(f)| < \infty$.
Most postcritically finite rational maps are {\em rigid}, in the sense defined in \S\ref{sec:Belyi};
 for example, $f$ is rigid if $|P(f)| > 4$.
When $f$ is rigid, it is conformally conjugate to a rational map $g$ defined over the algebraic closure of $\ratls$, in which case we have $P(g)\subset \proj^1(\ratlsbar)$.  In this paper we prove 
the converse:

\begin{theorem} \label{thm:Belyi1}
For any finite set $X\subset \proj^1(\ratlsbar)$ with $|X| \ge 2$, there exists a rigid, 
postcritically finite rational map such that 
$P(f)=X$.
\end{theorem}
We can also arrange that $P(f) \subset C(f)$, which implies that $f$ is hyperbolic.
The proof (\S\ref{sec:Belyi}) uses Belyi's theorem; consequently, the degree of the map $f$ 
we construct may be enormous, even when the set $X$ is rather simple.

In light of Theorem \ref{thm:Belyi1}, we formulate the following more precise question 
about dynamics on the postcritical set.

\begin{quest}
\label{question1} Let $X\subset \proj^1(\ratlsbar)$ be a finite set. Is every map $F:X\to X$ 
realized by a rigid rational map $f:P(f)\to P(f)$ with $P(f) = X$? 
\end{quest}
It is easy to see the answer can be {\em no} when $|X|=2$, and it is {\em yes} when
$|X|=3$ (for a proof, see \S\ref{sec:exs}).
We do not know of any negative answer when $|X|=4$, but most cases
are simply open. For example, we do not know the answer for $X=\{0,1,4,\infty\}$ and $F(x)=x$. 

We can, however, show that Question \ref{question1} has a positive answer if we allow 
perturbations of  $X$ as a subset of $\proj^1(\cx)$.

\begin{theorem}
\label{thm:combinatorics}
Let $F:X \to X$ be an arbitrary map defined on a finite set 
$X \subset \proj^1(\cx)$ with $|X| \ge 3$.  
Then there exists a sequence of rigid postcritically 
finite rational maps $f_n$ 
such that $|P(f_n)| = |X|$, 
\begin{displaymath}
	P(f_n) \to X \AND f_n | P(f_n) \to F |X 
\end{displaymath}
as $n \to \infty$.
\end{theorem}
The proof (\S\ref{sec:iter}) uses iteration on Teichm\"uller space, as in the proof of Thurston's 
topological characterization of postcritically finite rational maps
\cite{Douady:Hubbard:Thurston}.

In \S\ref{sec:Hur}, we establish  the following new Hurwitz-type result, 
which may be of interest
in its own right:
\begin{quote}
{\em Any collection of partitions $\cP$ with $|\cP|\geq 3$ 
can be extended to the passport $\cQ$ of a rational map, with $|\cQ| = |\cP|$}.
%{\em Any collection of partitions $\cP$ with $|\cP| \ge 3$
%can be extended to the passport of a rational map.}
\end{quote}
See Theorem \ref{thm:passport}.   (Here the {\em passport} of a rational map $f$ is 
the collection of partitions of $\deg(f)$ arising from the fibers over its critical values.)
We use this result to strengthen Theorem \ref{thm:combinatorics} in \S\ref{sec:local_degreez} 
by showing that we can 
also specify the multiplicity of $F$ at each point of $X$.

Other constructions of rational maps with specified postcritical 
sets are presented in \S\ref{sec:exs}.

\section{Every finite set of algebraic numbers is a postcritical set}
\label{sec:Belyi}

In this section we prove that any
finite set $X \subset \proj^1(\ratlsbar)$ with $2 \le |X| < \infty$ arises
as the postcritical set of a rigid rational map (Theorem \ref{thm:Belyi1}).
The result for rational maps follows easily from the following variant for 
polynomials.

\begin{theorem}
\label{thm:polyP}
For any finite set of algebraic numbers $X \subset \cx$
with $|X| \ge 1$, there exists a polynomial $f$ such that $P(f) \cap \cx = X$.
\end{theorem}

\bold{Polynomials.}
To begin the proof, we remark that it is easy to construct polynomials
with prescribed critical values.  It is convenient, when discussing a
polynomial $f$, to omit the point and infinity and let
$C_0(f) = C(f) \cap \cx$, and similarly for $V_0(f)$ and $P_0(f)$.

\begin{lemma}
\label{lem:cvs}
For any finite set $X \subset \cx$, there exists a polynomial $g$ with 
$V_0(g)=X$.
\end{lemma} 

\bold{Proof.} 
There are many ways to prove this result; for example, by induction on $|X|$, using the fact
that $V_0(f \compos g) = V_0(f) \cup f(V_0(g))$ and $V_0(z^2+a) = \{a\}$.  
\qed
\newline\noindent
(For a more precise result, see Corollary \ref{cor:polypass}.)

Let us say $\beta$ is a {\em Belyi polynomial} if
$V_0(\beta) \subset \{0,1\}$.
We will also use following result from \cite{Belyi:triple}: 

\begin{theorem}
\label{thm:Belyi}
For any finite set $X \subset \ratlsbar$, there exists a Belyi polynomial
such that $\beta(X) \subset \{0,1\}$.
\end{theorem}

\bold{Proof of Theorem \ref{thm:polyP}.} 
Our aim is to construct a polynomial $f$ such that $P_0(f)=X$.
This is easy if $|X| \le 1$, so assume $|X| \ge 2$.
Using Lemma \ref{lem:cvs}, choose a polynomial  $g$ such that
$V_0(g) = X$, and precompose with an
affine transformation so that $\{0,1\} \subset C_0(g)$.
Let $\beta$ be a Belyi polynomial such that $\beta(X) \subset \{0,1\}$.
Finally, let $f = g \circ \beta$.

We claim that $V_0(f) = X$.  Indeed, we have
\begin{displaymath}
	V_0(f) = V_0(g \circ \beta)  = V_0(g) \cup g(V_0(\beta)) = X \cup g(V_0(\beta)) ;
\end{displaymath}
but $V_0(\beta) \subset \{0,1\} \subset C_0(g)$, and therefore
$g(V_0(\beta)) \subset X$, so $V_0(f)=X$.
In particular $X \subset P_0(f)$.
But in fact $P_0(f) = X$, since
\begin{displaymath}
	f(X) = g \circ \beta(X)  \subset g(\{0,1\}) \subset X .
\end{displaymath}
\qed 

\bold{Hyperbolicity.}
Note that
\begin{displaymath}
	P_0(f) = X \subset \beta^{-1}(C_0(g)) \subset C_0(g \compos \beta) = C_0(f)
\end{displaymath}
in the construction above.
Thus every periodic point in $P(f)$ is superattracting,
and hence $f$ is hyperbolic.

\bold{Rigidity.} 
To deduce Theorem \ref{thm:Belyi1}, we must first
briefly discuss rigidity.
A postcritically finite map $f:\proj^1\to \proj^1$ is {\em rigid} if any postcritically finite
rational map $g$ uniformly 
close enough to $f$, and with $|P(g)|=|P(f)|$, is in fact 
conformally conjugate to $f$.
For any fixed $d$ and $n$, the rigid maps with $\deg(f)=d$ and $|P(f)|=n$
fall into finitely many conjugacy classes.

By a theorem of Thurston \cite{Douady:Hubbard:Thurston}, the only postcritically finite rational maps
$f$ that are {\em not} rigid are the flexible Latt\`es examples, which arise from the
addition law on an elliptic curve \cite{Milnor:Lattes}.  These flexible maps
have $|P(f)|=4$ and Julia set $J(f)=\proj^1$.  Consequently, 
any postcritically finite rational map with a 
periodic critical point (such as a polynomial) is rigid.

\bold{Proof of Theorem \ref{thm:Belyi1}.} 
Let $X\subset \proj^1(\ratlsbar)$ be a finite set with $|X| \ge 2$.
After a change of coordinates defined over $\ratlsbar$, we can assume that
$\infty \mem X$.  Then by Theorem \ref{thm:polyP},
there exists a polynomial $f$ with $P(f)=X$; and as we have just
observed, any postcritically finite polynomial is rigid.
\qed

\bold{Belyi degree and postcritical degree.}
Given a finite set $X \subset \ratlsbar$,
let $B(X)$ denote the minimum of $\deg(\beta)$ over all
Belyi polynomials with $\beta(X) \subset \{0,1\}$;
and let $D(X)$ denote the minimum of $\deg(f)$ over
all polynomials with $P_0(f) = X$.

Little is known about the general behavior of 
these degree functions
(in particular, lower bounds seem hard to come by);
however, the proof of Theorem \ref{thm:polyP}
in concert with Corollary \ref{cor:polypass} gives the relation:
\begin{displaymath}
	D(X) \le B(X) + |X| + 1.
\end{displaymath}
Both degree functions seem to merit further study.

\section{Contraction on Teichm\"uller space}
\label{sec:density}

The rational maps $f$ constructed in the proof of Theorem \ref{thm:Belyi1} all satisfy
$|f(P(f))|\le 3$.   
We next address the problem of realizing more general dynamics on $P(f)$.
Our construction will use iteration on Teichm\"uller space as in \cite{Douady:Hubbard:Thurston}.
This section gives the needed background;
for more details see \cite{Buff:Cui:Tan:survey}, \cite{Hubbard:book:T2}.

\bold{Teichm\"uller spaces.}
Given a finite set $A \subset \proj^1$ with $|A|=n$,
we let $\cT_A \isom \cT_{0,n}$ denote the
{\em Teichm\"uller space} of genus zero Riemann surfaces {\em marked} by $(\proj^1,A)$.

A point in $\cT_A$ is specified by another pair $(\proj^1,A')$ together with 
an orientation--preserving {\em marking} homeomorphism:
\begin{displaymath}
	\phi : (\proj^1, A) \arrow (\proj^1,A').
\end{displaymath}
Since $A'=\phi(A)$, the marking $\phi$ alone determines a point $[\phi] \mem \cT_A$.
Two markings $\phi_1, \phi_2$ determine the same point iff we can write $\phi_2 = \alpha \compos \phi_1 \compos \psi$,
where $\alpha \mem \Aut(\proj^1)$ and $\psi$ is isotopic to the identity {\em rel} $A$.

The cotangent space to $\cT_A$ at $(\proj^1,A')$ is naturally identified with the vector space $Q(\proj^1-A')$
consisting of meromorphic differentials $q = q(z) \,dz^2$ on $\proj^1$
with at worst simple poles on $A'$ and elsewhere holomorphic.
The Teichm\"uller metric corresponds to the norm
\begin{displaymath}
	\|q\| = \int_{\proj^1} |q| 
\end{displaymath}
on the cotangent space.

\bold{Pullback.}
Now let $F:\proj^1\to \proj^1$ be a {\em smooth} branched covering map with $\deg(F) \ge 2$.
The sets $C(F)$, $V(F)$ and $P(F)$ are defined just as for a rational map.

Consider a pair finite sets $A$ and $B$ in $\proj^1$ such that
\begin{displaymath}
	F(A) \cup V(F) \subset B .
\end{displaymath}
We then have a map of pairs
\begin{displaymath}
	F : (\proj^1,A) \arrow (\proj^1,B)	
\end{displaymath}
that is branched only over $B$.  By pullback of complex structures, we then obtain a holomorphic map
\begin{displaymath}
	\sigma_F : \cT_B \arrow \cT_A .
\end{displaymath}

\bold{Contraction.}
Let $\sigma : \cT_B \arrow \cT_A$ be a holomorphic map between Teichm\"uller spaces.
By the Schwarz lemma, $\|D\sigma\| \leq 1$ in the  Teichm\"uller metric.
We say $\sigma$ is {\em contractive} if $\|D\sigma\|<1$ at every point of $\cT_B$.  (We also say $\sigma$ is contractive if $|B|=3$, since
then its image is a single point.)
If $A=B$ and $\sigma$ is contractive, then $\sigma$ has at most one fixed point.

A branched covering $F$ is {\em contractive} if $\sigma_F$ is contractive.

The following result is well known and was a key step in the proof of Thurston's
rigidity theorem for postcritically finite rational maps; see
\cite[Prop. 3.3]{Douady:Hubbard:Thurston}.

\begin{prop}
\label{prop:con}
The map $F$ is contractive if and only if there is no
4-tuple $B_0 \subset B$ such that
\begin{equation}
\label{eq:con}
	F^{-1}(B_0) \subset A \cup C(F) .
\end{equation}
\end{prop}

\bold{Proof.}
We may assume $|B| \ge 4$.
Suppose $\|D\sigma_F\|=1$ at some point in $\cT_B$.
The coderivative of $\sigma_F$ at this point is given explicitly 
by a pushforward map of the form
\begin{displaymath}
	f_* : Q(\proj^1-A') \arrow Q(\proj^1-B') ,
\end{displaymath}
where $f$ is a rational map of the same topological type as $F$.
Since the domain of $f_*$ is finite dimensional,
there exists a nonzero $q \mem Q(\proj^1-A')$ such that $\|f_*q\| = \|q\|$.
The fact that there is no cancellation under pushforward
implies that $q$ is a locally a positive real multiple of
the pullback of $f_*q$.
In fact, since $\|f^*f_*q\| = \deg(f) \, \|f_*q\| = \|q\|$, we must have
\begin{displaymath}
	f^*f_* q = \deg(f) \, q.
\end{displaymath}
Now recall that any meromorphic quadratic differential on
$\proj^1$ has at least 4 poles.
Choose 4 points $B_0 \subset B$ such 
that $\phi(B_0) \subset B'$ is contained in the poles of $f_*q$.
Then the equation above implies that
the poles of $f^*f_*q$ lie in $A'$.
Thus any point in $F^{-1}(B_0)$ that does not lie in $A$
must be a critical point of $F$, giving condition
(\ref{eq:con}) above.

For the converse, suppose we have 4 points $B_0 \subset B$ satisfying (\ref{eq:con}).
Consider, at any point in $\cT_B$,
a quadratic differential $q$ with poles only at the 4 points marked by $B_0$.
Then $f^*q$ has poles only at points marked by $A$,
and hence it represents a cotangent vector to $\cT_A$.  
Since $\|f_*(f^*q)\| = \deg(f)\|q\| = \|f^*q\|$,
we have $\|D\sigma_F\|=1$ at every point in $\cT_B$.
\qed

%For the converse, assume there is a set of 4 points $B_0 \subset B$ satisfying \eqref{eq:con}.  Fix any point $\tau\in \cT_B$ represented by points $B' \subset \proj^1$, and let $f: (\proj^1, A') \to (\proj^1, B')$ be the holomorphic map with the topological type of $F$.  Let $v \in Q(\proj^1 - B')$ be a quadratic differential with exactly 4 simple poles at the points of $B_0'$, and set $q = f^*v \in Q(\proj^1 - A')$.   Then $f_*q = f_*f^*v = \deg(f) \, v$, so that $f^*f_*q = f^*(\deg(f)\, v)=  \deg(f)\, q$, and therefore $\|f_*q\| = \|q\|$.  Thus, $\|D\sigma_F\| = 1$ at $\tau$ and $F$ fails to be contractive.

\bold{Example:  dynamics.}
Let $f$ be a rational map with $|P(f)|=3$.
Note that a rational map is a special case of a smooth branched covering.
Consider a map of pairs
\begin{displaymath}
	f^k : (\proj^1,A) \arrow (\proj^1,B),
\end{displaymath}
such that $f^k(A) \cup P(f) \subset B$.

\begin{prop}
\label{prop:fk}
The pullback map 
\begin{displaymath}
	\sigma_{f^k} : \cT_B \arrow \cT_A
\end{displaymath}
is a contraction provided $\deg(f)^k > |A|$.
\end{prop}

\bold{Proof.}
We may assume $|B| \ge 4$.
Consider any 4-tuple $B_0 \subset B$.
Since $|P(f)| = 3$, we have a point $b \mem B_0 - V(f^k)$.
Then $f^{-k}(b)$ is disjoint from $C(f^k)$, and
$|f^{-k}(b)| = \deg(f)^k > |A|$, so we cannot have
$f^{-k}(B_0) \subset A \cup C(f^k)$.
\qed

\bold{Factorization.}
For later use, we record the following fact.
Suppose we have a factorization $F = F_1 \compos F_2$.
The pullback map can then be factored as
\begin{equation}
\label{eq:factor}
	\cT_B 
	\stackrel{\sigma_{F_1}}{\Arrow}
	\cT_C 
	\stackrel{\sigma_{F_2}}{\Arrow}
	\cT_A ,
\end{equation}
where $C = F_2(A) \cup V(F_2)$ .

\bold{Combinatorial equivalence.}
Finally we formulate the connection between fixed points
on Teichm\"uller space and rational maps, following Thurston.

Let $F$ and $G$ be a pair of postcritically finite
branched coverings of $\proj^1$.
An orientation--preserving homeomorphism of pairs 
\begin{displaymath}
	\phi : (\proj^1,P(F)) \arrow (\proj^1,P(G))
\end{displaymath}
gives a {\em combinatorial equivalence} between $F$ and $G$
if there is a second homeomorphism $\psi$, isotopic to
$\phi$ rel $P(F)$, making the diagram
\begin{displaymath}
\xymatrix{
 (\proj^1,P(F)) \ar[d]_{F} \ar[r]^{\psi} & (\proj^1,P(G)) \ar[d]^{G} \\
 (\proj^1,P(F))              \ar[r]^{\phi} & (\proj^1,P(G))  
}
\end{displaymath}
commute.  
Since $X=P(F)$ is forward invariant, $F$ determines a holomorphic map
\begin{displaymath}
	\sigma_F : \cT_X \arrow \cT_X.
\end{displaymath}
The following result follows readily from the
definitions (cf.  \cite[Prop. 2.3]{Douady:Hubbard:Thurston}):

\begin{prop}
\label{prop:fp}
A point $[Y] \mem \cT_X$ is fixed by $\sigma_F$ if and only
if there exists a rational map $f$ with $P(f)=Y$ such that
the marking homeomorphism
\begin{displaymath}
	\phi : (\proj^1,X) \arrow (\proj^1,Y)
\end{displaymath}
gives a combinatorial equivalence between $F$ and $f$.
\end{prop}

\section{Prescribed dynamics on $P(f)$}
\label{sec:iter}

In this section we prove Theorem \ref{thm:combinatorics}.
That is, given a finite set $X \subset \proj^1$ with $|X| \ge 3$, and a map $F : X \arrow X$, we will
construct a sequence of rigid rational maps $f_n$ such that $P(f_n) \arrow X$ and
\begin{displaymath}
	f_n|P(f_n) \arrow F|X .
\end{displaymath}
This means that for all $n \gg 0$, we can find homeomorphisms $\phi_n$ of $\proj^1$ such that
$\phi_n \arrow \id$, $\phi_n(X) = P(f_n)$, and $\phi_n$ conjugates $F|X$ to $f_n|P(f_n)$.

\bold{The setup.}
Let $h$ be a quadratic rational map with $J(h) = \proj^1$ and
\begin{displaymath}
	P(h) =\{0,1,\infty\}.
\end{displaymath}
(Explicitly, we can take $h(z) = (2/z - 1)^2$).
Note that $V(h^n) = P(h)$ for all $n \ge 2$.  
The map $h$ is expanding in the associated orbifold metric on $\proj^1$ 
(see e.g.\ \cite[Thm 19.6]{Milnor:book:dynamics},
	\cite[App. A]{McMullen:book:CDR}).

It is convenient to normalize so that $X$ contains $P(h)$.
Let $g$ be a polynomial fixing $0$ and $1$, such that
its finite critical values $V_0(g)$ coincide with $X-P(h)$.
(Such a polynomial exists by Lemma \ref{lem:cvs}).
Then
\begin{displaymath}
	V(g \compos h^n) = V(g) \cup g(V(h^n)) = X
\end{displaymath}
for all $n \ge 2$.
For later convenience, we also choose $g$ such that
\begin{equation}
\label{eq:g}
	\deg(g) = 3
	\;\;\;\text{if $|X|=4$} .
\end{equation}
(I.e.\ if $X = \{0,1,\infty,a\}$, we take $g$ to be
a cubic polynomial with $V_0(g)=\{a\}$.)

\bold{Approximation by branched covers.}
Since $J(h)=\proj^1$, the set $\bigcup_n h^{-n}(x)$
is dense for any $x \mem \proj^1$.   Using this fact, we can construct 
a sequence of homeomorphisms
\begin{displaymath}
	\phi_n : (\proj^1,X) \arrow (\proj^1,X_n),
\end{displaymath}
with $\phi_n \arrow \id$, such that
\begin{displaymath}
	F|X = g \compos h^n \compos \phi_n | X
\end{displaymath}
for all $n$. To do this we first pick, for each $x \mem X$, a nearby point $x'$ such that
$g \compos h^n(x') = F(x)$ and such that the map $x \mapsto x'$ is injective;
set $X_n = \{x' \st x \mem X\}$.  Then, we  choose a homeomorphism
$\phi_n$ close to the identity that moves $x$ to $x'$ for all $x \mem X$.
The larger $n$ is, the closer we can take $x'$ to $x$, and hence the closer
we can take $\phi_n$ to the identity.

\bold{Construction of rational maps.}
Next we observe that 
\begin{displaymath}
	F_n = g \compos h^n \compos \phi_n : (\proj^1,X) \arrow (\proj^1,X)
\end{displaymath}
is a smooth branched
covering map with $P(F_n) = X$.

Theorem \ref{thm:combinatorics} follows from:

\begin{theorem}
\label{thm:crux}
For all $n \gg 0$, $F_n$ is combinatorially equivalent to 
a rational map $f_n$; and suitably normalized, we have $f_n|P(f_n) \arrow F|X$.
\end{theorem}

\bold{Proof.}
Using the maps $\phi_n$, we can regard $X_n$ as points in $\cT_X$
such that $X_n \arrow X$.  
By construction, we have
\begin{displaymath}
	\sigma_{F_n}(X) = X_n ,
\end{displaymath}
and $d(X_n,X) \arrow 0$.  We wish to control the contraction of 
$\sigma_{F_n}$, and produce
a fixed point close to $X$.  That is, to complete the proof it suffices to show
there are point configurations $P_n \arrow X$ such that 
$\sigma_{F_n}(P_n)=P_n$.
For then, by Proposition \ref{prop:fp}, we have
a corresponding sequence of rational maps 
satisfy $P(f_n) = P_n$, and the marking homeomorphisms
transport $F|X$ to $f_n|P(f_n)$.

Choose $k$ such that $\deg(h)^k = 2^k > |X|+3$.
The crux of the matter is the factorization
\begin{displaymath}
	F_n = (g \compos h^k) \compos (h^{n-k} \compos \phi_n),
\end{displaymath}
valid for all $n \ge k$.  From this we obtain
a factorization of $\sigma_{F_n}$ as
\begin{displaymath}
        \cT_X
        \stackrel{\sigma_g}{\Arrow}
        \cT_B
        \stackrel{\sigma_{h^k}}{\Arrow}
        \cT_A
        \stackrel{\sigma_{h^{n-k} \circ \phi_n}}{\Arrow}
        \cT_X ;
\end{displaymath}
see equation (\ref{eq:factor}).
In this factorization, we have
\begin{displaymath}
	A = h^{n-k}(\phi_n(X)) \cup V(h^{n-k}).
\end{displaymath}
Since $|V(h^{n-k})| \le |P(h)| = 3$, we have
$\deg(h)^k = 2^k > 3 + |X| \ge |A|$, and hence $\sigma_{h^k}$ is a contraction
by Proposition \ref{prop:fk}.
Thus $\sigma_{g \compos h^k} = \sigma_{h^k} \compos \sigma_g$
is also a contraction. 
The amount of contraction at $P \mem \cT_X$ 
varies continuously with $P$.
Since the ball of radius $2$ about $X$ in $\cT_X$ is compact,
we can find a constant $\lambda$, independent of $n$, such that
\begin{displaymath}
	d(P,X) \le 2 \implies
	\|D\sigma_{F_n}(P)\|
	\le
	\|D\sigma_{g \compos h^k}(P)\|
	\le \lambda < 1.
\end{displaymath}

Let $X_n^i = \sigma_{F_n}^i(X)$, and let $\epsilon_n = d(X,X_n)$.
For all $n \gg 0$, we have $\epsilon_n \ll (1-\lambda)$.  Under this assumption, we can prove by induction that:
\begin{quote}
	(i) $d(X_n^i,X) \le 1$, and hence \\
	(ii) $d(X_n^i,X_n^{i-1}) \le \epsilon_n \lambda^i $.
\end{quote}
From this it follows that $X_n^i$ converges, as $i \arrow \infty$, to a fixed point $P_n$ of $\sigma_{F_n}$
with $d(X,P_n) \le \epsilon_n (1-\lambda)^{-1}$.  Since $\epsilon_n \arrow 0$, this completes the proof.
\qed

\bold{Proof of Theorem \ref{thm:combinatorics}.}
We just need to verify that $f_n$ is rigid.
But if $f_n$ is a flexible Latt\`es example,
then $|P(f_n)| = |X| = 4$ and hence $\deg(g) = 3$
by condition (\ref{eq:g}).
Moreover $\deg(f_n)$ is a square, contradicting the fact
that $\deg(F_n) = \deg(g) \deg(h)^n = 3 \cdot 2^n$.
\qed

\bold{Algorithmic solution.}
We have used the construction above as the basis for a practical
computer program that solves the approximation problem addressed
by Theorem \ref{thm:combinatorics}.  The iteration does not
quite take place on Teichm\"uller space; rather, we arrange that
$P_n$ is always close enough to $X$ that there is a unique
homeomorphism of $\proj^1$ close to the identity sending $P_n$ to $X$.

\section{Solution to a Hurwitz problem}
\label{sec:Hur}

In this section we establish the existence of polynomials and rational maps with constrained branch data. 
This will enable us to strengthen Theorem \ref{thm:combinatorics} by 
prescribing local degrees at points of $X$ as discussed in \S\ref{sec:local_degreez}. 

Our main result is:

\begin{theorem}
\label{thm:passport}
Let $\cP = (P_1,\ldots,P_n)$ be a finite list of partitions.
Then $\cP$ can be extended to the passport of a polynomial if $n \ge 2$,
and to the passport of a rational map if $n \ge 3$.
\end{theorem}

\bold{Partitions.}
A {\em partition} $P$ of $d \ge 0$ is a 
list of positive integers $(p_1,\ldots,p_s)$ such that $\sum p_i = d$.
The {\em trivial partition} has $p_i=1$ for $i=1,\ldots,d$.

Given a second partition $P' = (p_1',\ldots,p_t')$ of $d' = \sum p_i'$, we let
\begin{equation}
\label{eq:PQ}
	P+P' = (p_1,\ldots,p_s,p_1',\ldots,p_t')
\end{equation}
denote the combined partition of $d+d'$.

\bold{Passports.}
A {\em passport} of degree $d$ is a finite list $\cP = (P_1,\ldots,P_n)$ of 
nontrivial partitions of $d$.
For both partitions and passports, repetitions are allowed and the order in which elements appear is unimportant.  
We set $$c(\cP) = \sum_i (d - |P_i|).$$

\bold{Rational maps.}
For any $y$ in the target of a rational map $f$ we have a partition $P(f,y)$
 of $d=\deg(f)$ given by:
\begin{displaymath}
	\sum_{f(x) = y} \mult(f,x).
\end{displaymath}
The {\em passport of $f$} is the collection of partitions
\begin{displaymath}
	\cP(f) = (P(f,v_1),\ldots,P(f,v_n))
\end{displaymath}
arising from the critical values $\{v_1,\ldots,v_n\}$ of $f$.
(The other points in the target of $f$ yield trivial partitions.)
The number of critical points of $f$ mapping 
to a given point $y$, counted with multiplicity, is $d - |P(f,y)|$.  Hence
	$$c(\cP(f)) = 2 d - 2.$$

\bold{Branched coverings.}
The passport of a smooth branched covering map $F : S^2 \arrow S^2$
is defined similarly.   As is well known, any passport that can be realized topologically can
be realized geometrically.  More precisely, we have:  

\begin{prop}
\label{prop:Thom}
Let $F : S^2 \arrow S^2$ be a branched covering with $V(F)=\{v_1,\ldots,v_n\}$,
and let $X = \{x_1,\ldots,x_n\} \subset \proj^1$.
Then there exists a rational map $f : \proj^1 \arrow \proj^1$ with $\deg(f)=\deg(F)$ and $V(f)=X$ such that
\begin{displaymath}
	P(f,x_i) = P(F,v_i) \;\;\text{for $i=1,\ldots,n$}.
\end{displaymath}
In particular, $\cP(f)=\cP(F)$.
\end{prop}

\bold{Proof.}
Choose an orientation--preserving diffeomorphism $\phi: S^2 \arrow \proj^1$ such that $\phi(v_i) = x_i$ for $i=1,\ldots,n$.
Pulling the complex structure on $\proj^1$ back to $S^2$ via $\phi \compos F$, and applying the uniformization theorem,
we obtain a homeomorphism $\psi : S^2 \arrow \proj^1$ such that $f = \phi \compos F \compos \psi^{-1}$ is a holomorphic
branched covering, and hence a rational map (cf. \cite{Thom:ratl}).
\qed

\bold{Hurwitz problem.}
The {\em Hurwitz problem} is to characterize the
passports that arise from branched coverings of $S^2$.  A complete 
solution is not known; for background, see e.g. 
\cite[Ch. 5]{Lando:Zvonkin:book:graphs}.
Theorem \ref{thm:passport} addresses a variant of this problem where we allow the partitions to be {\em extended}. 

\bold{Extensions.}
A partition $Q$ {\em extends} $P$ if $Q=P+P'$ for some partition $P'$.
(For example, $1+3+5+7 = 16$ is an extension of $3+7=10$.)

Our main interest is in finite collections of partitions $\cP = (P_1,\ldots,P_n)$,
with repetitions allowed.
In this setting we say that $\cQ$ {\em extends} $\cP$ if, when suitably ordered,
we have $\cQ = (Q_1,\ldots,Q_n)$ and $Q_i$ extends $P_i$ for $i = 1,\ldots,n$.
In particular, if $\cQ$ extends $\cP$ then $|\cQ|=|\cP|$.

\bold{Polynomials.}
When $g : \cx \arrow \cx$ is a {\em polynomial}, its passport $\cP(g)$ 
is defined as the list of partitions $P(g,v_i)$ coming from the {\em finite} 
critical values of $g$.

The passports of polynomials are easily described.
In fact, by \cite[Prop. 5.2]{Edmonds:Kulkarni:Stong} we have:

\begin{theorem}
\label{thm:polypass}
A passport $\cP = (P_1,\ldots,P_n)$ of degree $d$ arises from a polynomial
$g$ if and only if
\begin{equation}
\label{eq:cd}
	c(\cP) = \sum_i (d - |P_i|) = d-1.
\end{equation}
\end{theorem}
The equation above is necessary because $g$ has $d-1$ critical points.
Applying Proposition \ref{prop:Thom}, we obtain:

\begin{cor}
\label{cor:polypass}
Let $X \subset \cx$ be a finite set such that $1 \le |X| < d$.
Then there exists a polynomial $g$ of degree $d$ whose critical values coincide with $X$.
\end{cor}

\begin{cor}\label{cor:poly}
Let $\cP = (P_1,\ldots,P_n)$ be a collection of partitions
with $n \ge 2$.  Then $\cP$ can be extended to the passport of a polynomial
of degree $d$ for all $d$ sufficiently large.
\end{cor}

\bold{Proof.}
It will be convenient to use exponential notation for repeated integers
(so $(1^d)$ is a partition of $d$).  

First, extend the partitions in $\cP$ so they are all nontrivial partitions of
the same integer $d$. Then $\cP$ is a passport.   If we extend $P_i$ to 
$P_i+(1)$ for all $i$, then $d$ increases by $1$ but $c(\cP)$ remains the same.  
Thus after a further extension of $\cP$, we can assume that $d-1 \ge c(\cP)$.  
If equality holds, we are done.

Otherwise, extend $P_i$ to $P_i+(3)$ for $i=1,2$, and to
$P_i+(1,1,1)$ for $i \ge 3$.  Then $d$ increases by $3$ while $c(\cP)$ increases by $4$.
By repeating this type of extension
until equality holds in equation (\ref{eq:cd}),
we obtain an extension of $\cP$ that arises from a polynomial of degree $d$.

Finally, suppose $\cP$ is the passport of a polynomial of degree $d$.
To complete the proof, we will show that
$\cP$ can be extended to a polynomial passport of degree $d+k$ 
for any $k \ge 2$.
To see this, just extend $P_1$ to $P_1+(k)$, $P_2$ to $P_2+(2,1^{k-2})$,
and $P_i$ to $P_i + (1^k)$ for $i \ge 3$.
\qed

\begin{theorem}\label{thm:ratl}
Any collection of partitions $\cP = (P_1,\ldots,P_n)$ 
with $n \ge 3$ can be extended to the passport of a rational map.
\end{theorem}

\bold{Proof.}
We divide the proof into two cases.

\bold{Case I.}  Assume $n \ge 4$.
Let $\cP_1 = (P_1,P_2)$ and let $\cP_2 = (P_3,\ldots,P_n)$.
By Corollary \ref{cor:poly}, we may assume that $\cP_1$ and $\cP_2$ are 
passports of polynomials $g_1$ and $g_2$ of the same
degree $d$.  

The complex plane can be naturally completed
to a closed disk
\begin{displaymath}
	D \isom \cx\cup S^1
\end{displaymath}
by adding a circle to represent the rays in $T_\infty \chat$.
Then each polynomial $g_i$ extends continuously to a proper map 
$D_i$ of degree $d$ on $D$,
satisfying $D_i(x) = dx$ for all $x \mem S^1 \isom \reals/\zed$.  

Now construct a branched covering $F : S^2\arrow S^2$
by gluing together two copies of $D$ to obtain a sphere,
and then setting $F =D_1$ on the first copy and $F=D_2$ on the second.
Then $\cP(F) = \cP$ by construction, and $F$ can be replaced by
a rational map by Proposition \ref{prop:Thom},
and the proof in this case is complete.

\begin{figure}[ht] 
   \centering
   \includegraphics[width=4.0in]{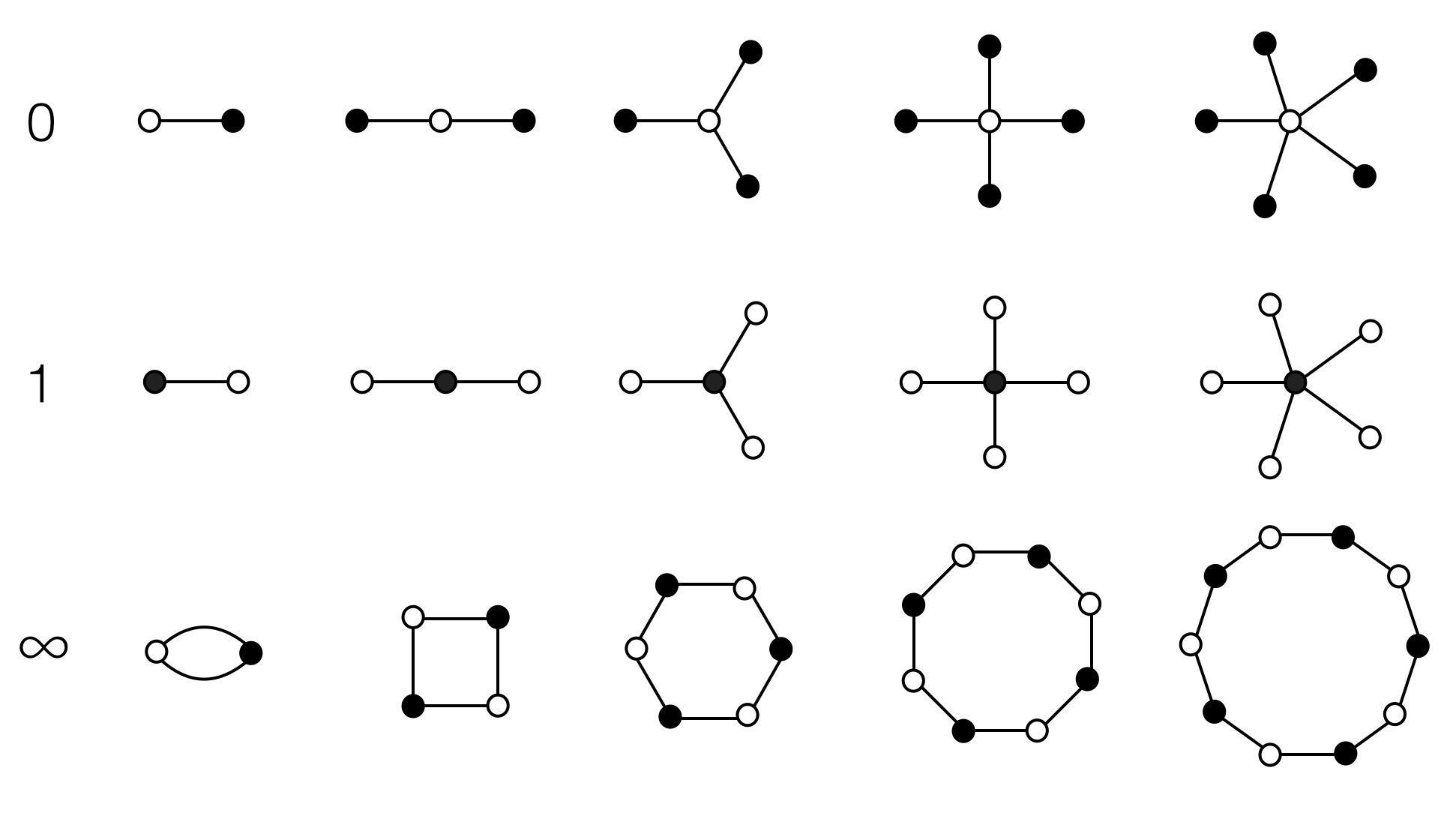} 
   \caption{\small Model dessins $D_0(m), D_1(m)$, and $D_\infty(m)$ 
   are drawn in rows for $1\leq m\leq 5$. }
   \label{dessinmodels}
\end{figure}

\begin{figure}[ht] 
   \centering
   \includegraphics[width=3in]{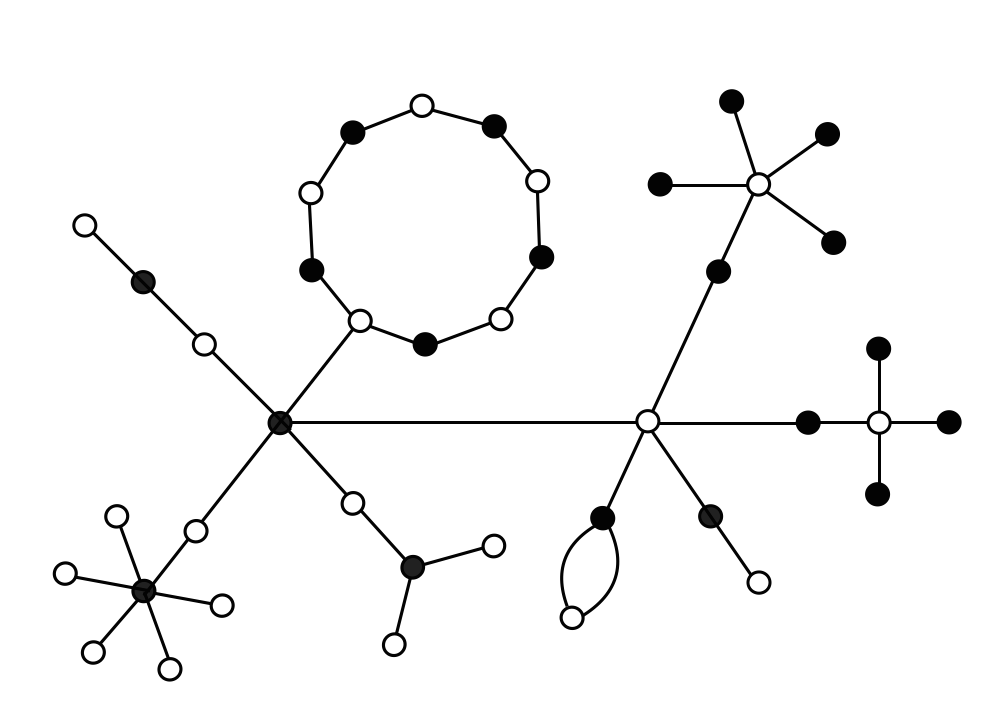} 
   \caption{\small A dessin d'enfant $D$ obtained by connecting 
   the components of $G$ to a single edge $[b,w]$, 
   which is drawn in the center. All four components 
   of $G_0$ are connected to $w$; all four components 
   of $G_1$ are connected to $b$; and one component of 
   $G_\infty$ is connected to $w$ 
   while the other is connected to $b$. 
   The corresponding rational map $f$ has degree $42$. } 
   \label{fancydancers}
\end{figure}

\bold{Case II.}  Now assume that $n=|\cP| = 3$. We will construct a rational map with $V(f)=\{0,1,\infty\}$ 
such that $\cP(f)$ extends $\cP$.

For convenience, let us index the elements $P_i$ of $\cP$ 
by $i \mem I = \{0,1,\infty\}$.  
After replacing $P_i$ by the extension $P_i+(2)$, if necessary,
we may assume that every partition in $\cP$ is nontrivial.
 
To construct $f$, it suffices to give the topological data of a
{\em dessin d'enfant} $D \subset \cx$ 
\cite{Schneps:book:dessins}, \cite{Lando:Zvonkin:book:graphs}. 
A dessin is connected graph with vertices of two colors, embedded in the plane,
arising as the preimage of the interval $[0,1]$ 
under a branched covering $F:\proj^1\to\proj^1$ with $V(F) \subset \{0,1,\infty\}$. 
We adopt the convention that
the white vertices of $D$ map to $0$ and the black vertices 
map to $1$.  The components of $\proj^1-D$ are called the {\em faces} of $D$,
and each face contains a unique point $z$ such that $F(z)=\infty$.

Consider, for each $m \ge 1$, the three types of model dessins
\[
	D_0(m), \;\;\; D_1(m),\AND D_\infty(m)
\]
shown in Figure \ref{dessinmodels}.  
(These graphs correspond to the rational maps
$f_0(z) = z^m$, $f_1(z) = z^m+1$ and $f_\infty(z) = (z^m + z^{-m} +2)/4$.)
The barycenter $z$ of $D_i(m)$ is a vertex for $i=0$ or $i=1$, and it lies in a face if $i=\infty$.  
In all three cases, $F(z)=i$ and $\mult(F,z)=m$.  

Let $P_0  = (a_1,\ldots,a_n)$, and let $G_0 \subset \cx$ be a planar graph with
$n$ components, such that its $j$th component is isomorphic to $D_0(a_j)$.
Such a graph is unique up to planar isotopy.
Construct $G_1 \subset \cx$ and $G_\infty \subset \cx$ similarly,
using the partitions $P_1,P_\infty$ and the models $D_1(m), D_\infty(m)$.
When constructing $G_\infty$, we take care not to nest two components of type $D_\infty(m)$.
It is easy to arrange that the graphs $G_i$ are contained in disjoint disks in the plane.
Let $G = G_0 \cup G_1 \cup G_2$.  Then every vertex of $G$ is incident to the
unique unbounded component $U$ of $\cx - G$.

To obtain $D$, we add new edges between vertices of opposite colors
to make $G$ connected.
In the process we take care not to alter the valence at the center of each component of $G$.
This can be done in many ways. 

For example, begin by introducing a new edge $[b,w] \subset U$, with one 
white vertex $w$ and one black vertex $b$. 
Then connect each component of $G_0$ to $w$, connect each component of
$G_1$ to $b$, and connect each component of $G_\infty$ to either 
$w$ or $b$ as in Figure \ref{fancydancers}, making sure that these new edges 
do not cross. 

The resulting connected graph $D \subset \cx$ is then a dessin d'enfant for 
a rational map $f$ with $V(f) = \{0,1,\infty\}$.
By the construction of the graph $G_i$, the partition $\cP(f,i)$
extends $P_i$ for each $i \mem V(f)$, and hence $\cP(f)$ extends $\cP$.  
\qed

\bold{Proof of Theorem \ref{thm:passport}.} 
The statements for polynomials and for rational maps are covered by
Corollary \ref{cor:poly} and Theorem \ref{thm:ratl} respectively.
\qed

\bold{Remarks and references.}
The proof of Theorem \ref{thm:ratl} for $n \ge 4$
is based on the fact that the passports for a pair of
polynomials of degree $d$ can be combined to
give the passport of a rational map.
This result also appears in
\cite[Remark on p. 785]{Edmonds:Kulkarni:Stong}
and \cite[Prop. 10]{Baranski:hurwitz}.
The branched covering built from a pair of polynomials is
called their {\em formal mating} in \cite{Tan:matings}. 

\section{Prescribed critical points in $P(f)$}
\label{sec:local_degreez}

In this section we strengthen Theorem \ref{thm:combinatorics} 
by showing that we can construct a rational map
with prescribed critical points in $P(f) \approx X$.
We will also give a similar result for polynomials.

\bold{Multiplicities.}
To make the first statement precise, recall that
$f$ has a critical point of order $\mult(f,x)-1$
at each point $x \mem \proj^1$.  
We can regard the multiplicity as a map
$\mult(f) : \proj^1 \arrow \zed_+ = \{1,2,3,\ldots\}$.

Our aim is to show:

\begin{theorem}
\label{thm:mult}
 Let $F:X \to X$ and $M:X\to \zed_+$ be 
 arbitrary maps defined on a finite set $X \subset \proj^1$ 
 with $|X| \ge 3$.  Then there exists a sequence of 
 rigid, postcritically finite rational maps $f_n$ 
 such that $|P(f_n)| = |X|$, 
\begin{displaymath}
	P(f_n) \arrow X, \;\;\; f_n | P(f_n) \to F|X, 
	\AND \mult(f_n) | P(f_n) \to M|X
\end{displaymath}
as $n \to \infty$.
\end{theorem}

\bold{Proof.} 
The argument is a modification of the proof of Theorem \ref{thm:combinatorics}
given in \S\ref{sec:iter}. 

As in that section, the construction is based on a pair of rational maps $g$ and $h$.
Let $h$ be a quadratic rational map with $J(h)=\proj^1$ and $|P(h)|=3$.

To construct $g$, we first associate to each $x \mem X$ the partition
\begin{displaymath}
	P_x = (M(y_1),\ldots,M(y_n)) \;\;\;
\end{displaymath}
where $F^{-1}(x) = \{y_1,\ldots,y_n\} \subset X$.
Then $\sum_X |P_x| = |X|$.
(Note that if $x \notmem F(X)$, then $P_x$ is the empty partition and $|P_x|=0$.)

Let $P_x' = P_x + (1)$ be the partition obtained by padding $P_x$ with an extra $1$
at the end (notation as in equation (\ref{eq:PQ})), so $|P_x'| = |P_x|+1$.
By Proposition \ref{prop:Thom} and Theorem \ref{thm:ratl}, there is a rational map $g$ with critical values 
$V(g)=X$ such that $P(g,x)$ extends $P_x'$ for all $x \mem X$.  Since $|X| \ge 3$, we then have
\begin{equation}
\label{eq:gi}
	|g^{-1}(X)|  = \sum_X |P(g,x)| \ge \sum_X |P_x'| = 2|X| \ge |X|+3.
\end{equation}

Next, we construct an injective map
\begin{displaymath}
	\iota : X \arrow g^{-1}(X)  \subset \proj^1
\end{displaymath}
such that
\begin{equation}
\label{eq:iota}
	F(y) = g(\iota(y)) 
	\AND
	M(y) = \mult(g,\iota(y))
\end{equation}
for all $y \mem X$.  
To define $\iota(y)$, let $x = F(y)$ and recall that: 
$y$ determines a point $M(y) \mem \cP_x$, we have an inclusion $\cP_x \subset P(g,x)$, and 
 there is a bijection $P(g,x) \isom g^{-1}(x)$ which labels points by their multiplicities.
We define $\iota(y)$ to be the image of $M(y)$ under the composition $\cP_x \subset P(g,x) \isom g^{-1}(x)$.

Since $|P(h)|=3$ and equation (\ref{eq:gi}) holds, we can choose
$\alpha \mem \Aut(\proj^1)$ such that 
$\alpha(P(h)) \subset g^{-1}(X) - \iota(X)$.
Upon replacing $g$ and $\iota$ with 
$g \compos \alpha$ and $\alpha^{-1} \compos \iota$,
equation (\ref{eq:iota}) continues to hold, and we then have
\begin{displaymath}
	P(h) \subset g^{-1}(X) - \iota(X) .
\end{displaymath}

Now recall that $J(h) = \proj^1$ and hence
\begin{equation}
\label{eq:XJ}
	X \subset J(h) .
\end{equation}
Consequently, for any $x \mem X$, the inverse orbit of $\iota(x)$ under 
$h$ accumulates on every point of $X$.
Moreover, if $h^n(z) = \iota(x)$,
then $\mult(h^n,z) = 1$, since $\iota(x) \notmem P(h)$.
Thus we can find a sequence of injective maps $\phi_n : X \arrow X_n$,
converging to the identity, such that 
\begin{equation}
\label{eq:phin}
	\mult(h^n)|X_n = 1 \AND h^n(\phi_n(x)) = \iota(x).
\end{equation}
Extend $\brackets{\phi_n}$ to a sequence of homeomorphisms of $\proj^1$ converging to the identity as $n \arrow \infty$,
and let
\begin{eqnarray}\label{eqn:branched}
	F_n= g \compos h^n \compos \phi_n : (\proj^1,X) \longrightarrow (\proj^1,X).
\end{eqnarray}
Then by equation (\ref{eq:iota}), we have 
\begin{equation}
\label{eq:Fn}
	F_n|X = F
	\AND 
	\mult(F_n)|X = M|X .
\end{equation}

Finally, we reiterate the proof of Theorem \ref{thm:crux} to convert
the postcritically finite branched covers $F_n$ into rational maps $f_n$
with $P(f_n) \arrow X$.
Since $f_n$ and $F_n$ are conjugate on their postcritical sets, we then have
\begin{displaymath}
	f_n|P(f_n) \arrow F|X \AND \mult(f_n)|P(f_n) \arrow M|X
\end{displaymath}
by equation (\ref{eq:Fn}).
We can also ensure, by our choice of $g$, that $\deg(f_n)$ is not a square, and hence
the postcritically finite maps $f_n$ are rigid for all $n \gg 0$. 
\qed

\bold{The polynomial case.}
We conclude by presenting a variation of Theorem \ref{thm:mult} for polynomials.

\begin{theorem}\label{thm:mult:poly}
Let $F:X \to X$ and $M:X\to \zed_+$  be 
 arbitrary maps defined on a finite set $X \subset\cx$ 
 with $|X| \ge 2$.  Then there exists a sequence of 
postcritically finite polynomials  $f_n$ 
 such that $|P_0(f_n)| = |X|$, 
\begin{displaymath}
	P_0(f_n) \arrow X, \;\;\; f_n | P_0(f_n) \to F|X, 
	\AND \mult(f_n) | P_0(f_n) \to M|X
\end{displaymath}
as $n \to \infty$.
\end{theorem}

\bold{Prescribed Julia sets.}
For the proof, we will need a polynomial $h$ whose Julia set contains $X$
(to play the role of the rational map $h$ with $J(h)=\proj^1$ in the proof of
Theorem \ref{thm:mult}).
It suffices to treat the case where $X \subset \ratlsbar$, since 
$\ratlsbar$ is dense in $\cx$. 

\begin{lemma}\label{lem:h:poly}
Given any finite set $X\subset \ratlsbar$, there is a polynomial $h$ so that 
\[
|P_0(h)|=2 \AND X \subset J(h).
\]
\end{lemma}

\bold{Proof.} 
By Theorem \ref{thm:Belyi}, there is a polynomial $\beta$ with 
$\beta(X)\cup V_0(\beta) \subset \{0,1\}$.  We can assume that 
$\deg(\beta) > 1$ and $V_0(\beta) = \{0,1\}$ (for example, by taking 
a Belyi polynomial for a larger set that contains $X$). 

There are two distinct points $a,b \mem \beta^{-1}(\{0,1\})$
that are not critical points of $\beta$. Indeed, the set $\beta^{-1}(\{0,1\})$ 
consists of $2d$ points, counted with multiplicity, and at 
most $2|C_0(f)|=2d-2$ of these are accounted for by critical points. 
        
Let $\alpha \mem \Aut(\cx)$ be an affine transformation sending
the ordered pair
$(a,b)$ to $(0,1)$, and set $h = \alpha \compos \beta$.
We then have:
\begin{displaymath}
	V_0(h)= \{a,b\},
	\;\;\; h(\{a,b\})\subset \{a,b\}, 
	\AND h(X)\subset \{a,b\}.
\end{displaymath}
The first two properties imply that $P_0(h) = \{a,b\}$, and in particular
$P_0(h)$ is disjoint from $C_0(h)$.  It follows that $P_0(h)$ is contained in the
Julia set of $h$.  Since the Julia
set is totally invariant and $h(X) \subset \{a,b\}$, we have 
$X \subset J(h)$ as well.
\qed

\bold{Proof of Theorem \ref{thm:mult:poly}.}
We may assume $X \subset \ratlsbar$.
Let $h$ be a polynomial associated to $X$ as in Lemma \ref{lem:h:poly}.
Let $(P_x' \st x \mem X)$ be the family of partitions constructed in
the proof of Theorem \ref{thm:mult}, and let
$g$ be a polynomial with $V_0(g)=X$,
provided by Proposition \ref{prop:Thom} and Corollary \ref{cor:poly},
such that $P(g,x)$ extends $P_x'$ for all $x \mem X$.
We can now simply repeat the proof of Theorem \ref{thm:mult},
using the mappings $g$ and $h$
to obtain the desired polynomials $f_n$.
\qed

\section{Alternative constructions}
\label{sec:exs}

In this section we discuss alternative constructions of rational maps
with postcritical sets satisfying $|P(f)| \le 4$.

\bold{The case $|P(f)|\leq 3$.}
The only rational maps with $|P(f)|=2$
are those that are conformally conjugate to $z\mapsto z^{\pm d}$. 
Question \ref{question1} then  
has a negative answer if $X\subset \proj^1$ has cardinality $2$; 
that is, $F:X\to X$ is realized by $f:P(f)\to P(f)$ if and only if $F$ is 
bijective. 

\begin{figure}[ht] 
   \centering
   \includegraphics[width=3in]{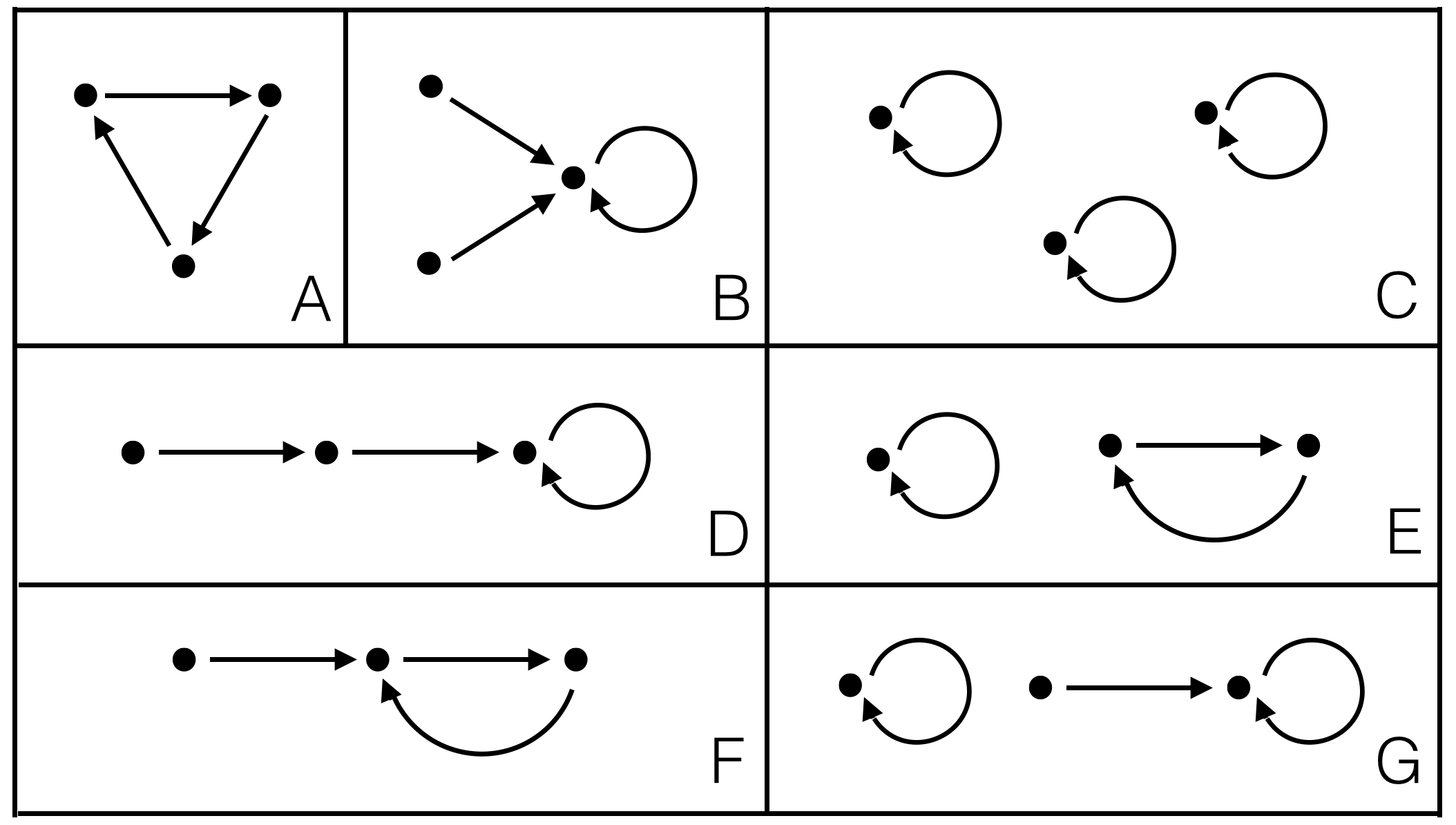} 
   \caption{\small This table contains all seven possible 
   graphs representing maps $X\to X$ when $|X|=3$.
 Each graph can be realized by a rational map $f:P(f)\to P(f)$.
Explicit examples are given by $A(z)= 1-1/z^2$; 
$B(z)= (z-\alpha)^3/(z-1+\alpha)^3$ where $\alpha^2-\alpha+1=0$;
$C(z)= z^2(3-2z)$; $D(z)=(1-2/z)^2$; $E(z)=z^2-1$;
 and $F(z)=(2z-1)^2/(4z(z-1))$; and $G(z)=z^2-2$.
In each case, the degree is minimal among all rational maps
realizing the given dynamics on $X$.
}
      \label{table}
\end{figure}

It is easy to see that Question \ref{question1} has a positive
answer whenever $|X|=3$.  
Indeed, we can first normalize so that
$X = \{0,1,\infty\}$.  Then, up to reordering the points of $X$,
there are only seven possibilities for $F$. 
A concrete rational map realizing
each one is given in Figure \ref{table}.

\bold{The case $|P(f)|=4$:  rigid Latt\`es maps.}
Let $Q \subset \half$ denote the set of $\tau$ in the 
upper halfplane that are quadratic over $\ratls$,
and let 
\begin{displaymath}
	\lambda : \half \arrow \half/\Gamma(2) 
	\isom \proj^1 - \{0,1,\infty\}
\end{displaymath}
be the universal covering map.  
Rigid Latt\`es maps are precisely those which are covered by
the action of complex multiplication on an elliptic curve $E$
(up to translation by a point of order in $E[2]$);
see \cite[Lemmas 4.3 and 4.4]{Milnor:Lattes}.
Every such elliptic curve has the form
$E = \cx/(\zed \oplus \zed \tau)$ with $\tau \mem Q$.
Using these rigid maps, one can
explicitly construct $f$ with $P(f) = \{0,1,\infty,z\}$
for all $z$ in the dense set $\lambda(Q) \subset \proj^1$.
On the other hand, 
\begin{eqnarray}\label{eq:L}
	L = \lambda(Q)
\end{eqnarray}
is very small subset of $\ratlsbar$.  
For example, $L \cap \zed = \{-1,2\}$, since the $j$--invariant
\begin{displaymath}
	j = \frac{256 (1-\lambda+\lambda^2)^3}{\lambda^2 (1-\lambda)^2}
\end{displaymath}
is an algebraic integer for all $\lambda \mem L$
(see e.g.\ \cite[\S II.6]{Silverman:book:advanced} for details).
So the postcritical sets arising from rigid Latt\`es examples
are insufficient to complete the proof of 
Theorem \ref{thm:Belyi1} in the case $|P(f)|=4$.

\bold{The case $|P(f)|=4$:  dynamics on moduli space.}
Our second construction of rational maps with given postcritical sets is 
based on \cite{Koch:pcf}.  There, the author builds maps on 
moduli space $g:\mathcal M_{0,n}\dashrightarrow \mathcal M_{0,n}$ 
whose periodic points correspond to postcritically finite rational 
maps on the Riemann sphere.  

For $n=4$, the moduli space $\mathcal M_{0,4}$ can 
be identified with $\proj^1-\{0,1,\infty\}$, and 
$g(z)=(1-2/z)^2 $ is an example of one such map on moduli space.  
Each periodic point $x$ of period $m$ gives rise 
to a postcritically finite polynomial of degree $2^m$ 
 with postcritical set $\{0,1,\infty, x\}$.  
 Let $K$ denote the set of all periodic points of $g$ in 
 $\proj^1 - \{0,1,\infty\}$.  The set $K$ 
lies in $\ratlsbar$, and it is dense in the Julia set $J(g)$. 
Moreover, we have $J(g)=\proj^1$ since $g$ is a Latt\`es map, so 
$K\subset \proj^1$ provides a dense set of postcritical sets 
arising from polynomials. 

Interestingly, the intersection $L\cap K$ is finite
(where $L$ is defined by equation (\ref{eq:L})). 
Indeed, $K$ is a set of bounded Weil height, as a consequence 
of the existence of a canonical height for $g$ \cite{Call:Silverman:heights}. 
But the $j$-invariants of elliptic curves with complex multiplication 
have no infinite subsets of bounded height \cite[Lemma 3]{Poonen:hecke}, 
so neither does $L$.

\bibliographystyle{math}
\bibliography{math}

\end{document}